# FEM-PIKFNNs for underwater acoustic propagation induced by structural vibrations in different ocean environments

Qiang Xi[b],   Zhuojia Fu[a,c,1],   Wenzhi Xu[c],   Mi-An Xue[a,b],   Jinhai Zheng[a,b]

[a]Key Laboratory of Ministry of Education for Coastal Disaster and Protection, Hohai University, Nanjing 210098, China

[b]College of Harbour, Coastal and Offshore Engineering, Hohai University, Nanjing 210098, China

[c]Center for Numerical Simulation Software in Engineering and Sciences, College of Mechanics and Materials, Hohai University, Nanjing 211100, China

**Abstract:** In this paper, a novel hybrid method based on the finite element method (FEM) and physics-informed kernel function neural networks (PIKFNNs) is proposed and applied to the prediction of underwater acoustic propagation induced by structural vibrations in the unbounded ocean, deep ocean and shallow ocean. In the hybrid method, PIKFNNs are a class of improved shallow physics-informed neural networks (PINNs) that replace the activation functions in PINNs with the physics-informed kernel functions (PIKFs), thereby integrating prior physical information into the neural network model. Moreover, this neural network circumvents the step of embedding the governing equations into the loss function in PINNs, and requires only training on boundary data. By using the Green's functions as the PIKFs and the structural-acoustic coupling response information obtained from the FEM as boundary training data, the PIKFNNs can inherently capture the Sommerfeld radiation condition at infinity, which is naturally suitable for predicting ocean acoustic propagation. Numerical experiments demonstrate the accuracy and feasibility of the FEM-PIKFNNs in comparison with the true solutions and finite element results.

---
[1]Corresponding authors: paul212063@hhu.edu.cn

## 1. Introduction

As the most popular numerical algorithm in the fields of science and engineering, the finite element method (FEM) [1-4] is widely applied to structural-acoustic coupling analysis in finite domains. However, the FEM will encounter a computational bottleneck where the computational cost increases rapidly as the domain size increases in the infinite domain ocean acoustic propagation simulation. In contrast, PIKFs-based numerical methods [5-7], such as boundary element method (BEM) [8,9], singular boundary method (SBM) [10-12] and method of fundamental solutions (MFS) [13-16], are preferred in the ocean acoustic simulation due to the employed PIKFs (Green's functions) inherently satisfy the governing equation of the ocean acoustic propagation and the Sommerfeld radiation condition at infinity. Moreover, these PIKFs-based numerical methods can also effectively avoid the truncation treatment required by the FEM for infinite domains. To take advantage of the individual strengths of the aforementioned methods, an increasing number of researchers have devoted their efforts to the study of coupling methods that integrate the FEM with these PIKFs-based numerical methods, and have proposed various coupling methods such as the FEM-BEM [17,18], FEM-SBM [19,20], and FEM-MFS [21,22] for vibration [23] and sound propagation analysis in infinite domains.

Over the past two decades, neural networks [24-29] have emerged as powerful machine learning tools with remarkable potential in the field of scientific computing. By employing multi-layer network connections and weight adjustments, they can capture intricate relationships among data and offer a fresh perspective for solving partial differential equations (PDEs). Among them, the physics-informed neural networks (PINNs) [30-34] and deep Galerkin method (DGM) [35,36] are two typical neural networks. They incorporate the information embedded in the governing equations into neural networks by using appropriate loss functions and training strategies. This ensures that the trained results of the neural network model not only fit the existing data but also satisfy the governing equations. These approaches have shown remarkable potential for solving high-dimensional, complex geometric and nonlinear problems. However, the iterative solving of multi-layer network structures during the training process has resulted in neural networks requiring more computational resources than traditional numerical methods. To improve the computational efficiency, the shallow PINNs [37,38] is proposed. It requires only iterative solutions for a single hidden layer network structure during each training step, which can significantly reduce the computational burden of the neural network with the same number of neurons. It is worth noting that the early proposed radial basis function neural networks (RBFNNs) [39,40] can also be considered as a type of shallow PINNs using the radial basis function (RBF) [41] as the activation function.

Recently, inspired by the fundamental ideas of PIKFs-based numerical methods and shallow PINNs, physics-informed kernel function neural networks (PIKFNNs) [42] have been proposed. PIKFNNs are a class of improved shallow PINNs that utilize PIKFs as neuron functions, thereby naturally incorporating prior physical information into the neural network model. Furthermore, PIKFNNs benefit from the inherent satisfaction of the governing equations by PIKFs. There is no need to explicitly embed the governing equations in the loss function. It only requires training to fit the boundary data, which significantly reduces the amount of training data, thereby achieving a more efficient training process. It is worth noting that when Green's function is used as the PIKFs, the PIKFNNs can automatically capture the Sommerfeld radiation condition at infinity. This property makes PIKFNNs promising for applications in the ocean acoustic propagation simulation.

Based on the above statements, this paper proposes the hybrid FEM-PIKFNNs and applies it to the prediction of underwater acoustic propagation induced by structural vibrations in the unbounded ocean, deep ocean and shallow ocean. The brief outline of this paper is given as follows. The theoretical and numerical model of the hybrid FEM-PIKFNNs for underwater acoustic propagation induced by structural vibrations in different ocean environments is provided in Section 2. In Section 3, the feasibility and efficacy of the hybrid FEM-PIKFNNs are presented under three numerical examples. The main findings are summarized in Section 4.

## 2. Theoretical and numerical model

In this section, the hybrid finite element method-physics-informed kernel function neural networks (FEM-PIKFNNs) is proposed to predict the underwater acoustic propagation induced by structural vibrations in the unbounded ocean, deep ocean and shallow ocean.

### 2.1. Finite element method

Based on the theory of linear elasticity, the finite element discrete equation for the time-harmonic vibration of an elastic structure is given by

$$\left(\mathbf{K}_a + i\omega\mathbf{C}_a - \omega^2\mathbf{M}_a\right)\mathbf{u} = \mathbf{F}_a, \tag{1}$$

where $\mathbf{K}_a$, $\mathbf{C}_a$ and $\mathbf{M}_a$ are the stiffness matrix, damping matrix and mass matrix of the elastic structure, $i = \sqrt{-1}$, $\omega$ is the angular frequency, $\mathbf{u}$ represents the displacement of the elastic structure, and $\mathbf{F}_a$ denotes the structural excitation force.

Similar to the structural vibration, the finite element discrete equation of the underwater acoustic propagation can be expressed as

$$\left(\mathbf{K}_b + i\omega\mathbf{C}_b - \omega^2\mathbf{M}_b\right)\mathbf{p} = \mathbf{F}_b, \tag{2}$$

where $\mathbf{K}_b$, $\mathbf{C}_b$ and $\mathbf{M}_b$ represent the stiffness matrix, damping matrix and mass matrix of the underwater acoustic propagation, $\mathbf{p}$ is the underwater sound pressure, and $\mathbf{F}_b$ is the acoustic excitation force.

In the ocean environment, the vibration of elastic structure and underwater acoustic wave must satisfy the conditions of force equilibrium and displacement continuity, the resulting structural-acoustic coupling matrix can be written in the following form

$$\begin{bmatrix} \mathbf{K}_a + i\omega\mathbf{C}_a - \omega^2\mathbf{M}_a & \mathbf{K}_c \\ -\omega^2\mathbf{M}_c & \mathbf{K}_b + i\omega\mathbf{C}_b - \omega^2\mathbf{M}_b \end{bmatrix} \begin{bmatrix} \mathbf{u} \\ \mathbf{p} \end{bmatrix} = \begin{bmatrix} \mathbf{F}_{ac} \\ \mathbf{F}_{bc} \end{bmatrix}, \tag{3}$$

where $\mathbf{F}_{ac}$ and $\mathbf{F}_{bc}$ are the coupling excitation forces of the elastic structure and underwater acoustic wave, $\mathbf{K}_c$ indicates the coupling stiffness matrix, $\mathbf{M}_c = -\rho_0 \mathbf{K}_c^T$, and $\oplus^T$ denotes the transpose of the matrix.

In the finite element implementation, ocean environments are truncated into finite domains by setting perfectly matched layers and ocean boundaries, as shown in Fig. 1. Then, the FEM mesh is generated in the finite domain, and the near-field underwater sound pressure induced by structural vibrations can be calculated by the FEM.

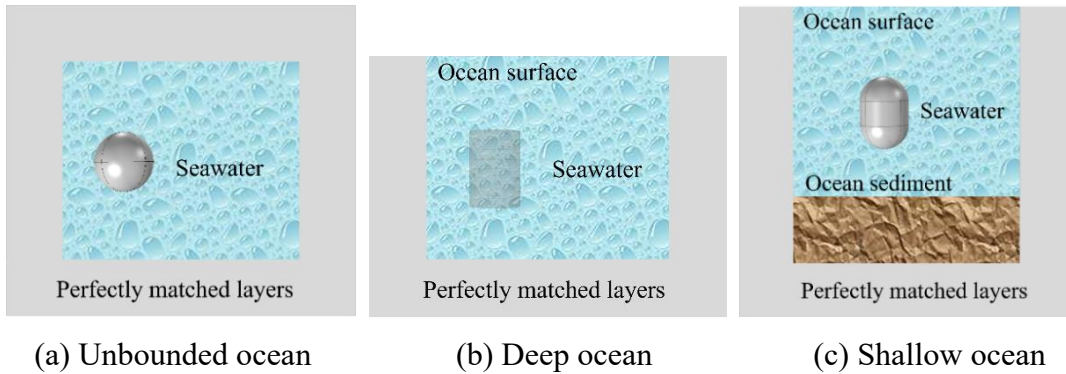

(a) Unbounded ocean      (b) Deep ocean      (c) Shallow ocean

Fig. 1. Sketch of the truncated finite domain.

## 2.2. Physics-informed kernel function neural networks

In this subsection, the near-field underwater sound pressure information calculated by the FEM is used as the training data of the PIKFNNs for predicting the underwater acoustic propagation in the unbounded ocean, deep ocean and shallow ocean.

### 2.2.1. Physics-informed kernel functions

Based on the linear acoustic wave theory, the governing equation for the sound pressure induced by the structural vibration in the ocean environment is given by

$$\Delta p(\mathbf{x}) + k^2 p(\mathbf{x}) = 0, \qquad (4)$$

in which $p$ denotes the sound pressure at the position $\mathbf{x}$, $\Delta$ indicates the Laplacian, and $k$ is the wave number. The sound pressure training sample is prescribed as

$$p(\mathbf{x}) = \bar{p}(\mathbf{x}), \qquad (5)$$

where $\bar{p}(\mathbf{x})$ represents the underwater sound pressure data calculated by the FEM. Moreover, the underwater sound pressure in all three ocean environments must satisfy the Sommerfeld radiation condition at infinity. It is worth noting that the underwater sound pressure in the deep and shallow oceans must also satisfy the additional reflection conditions imposed by the ocean surface or ocean sediment.

To construct the PIKFNNs, the PIKFs of underwater acoustic propagation in the unbounded ocean, deep ocean and shallow ocean are presented. In the unbounded ocean, the PIKF of underwater acoustic propagation is defined as

$$\Psi_1(\mathbf{x}, \mathbf{s}) = \frac{e^{ikR_a}}{R_a}, \qquad (6)$$

where $R_a = \sqrt{(x-s^x)^2 + (y-s^y)^2 + (z-s^z)^2}$ represents the distance between the sample point $\mathbf{x} = (x, y, z)$ and source point $\mathbf{s} = (s^x, s^y, s^z)$.

In the deep ocean, the propagation of radiated acoustic waves induced by structural vibrations is affected by the reflection effect from the ocean surface. Its PIKF can be represented as

$$\Psi_2(\mathbf{x}, \mathbf{s}) = \frac{e^{-ikR_a}}{R_a} - \frac{e^{-ikR_b}}{R_b}, \qquad (7)$$

where $R_b = \sqrt{(x-s^x)^2 + (y-s^y)^2 + (2h-z-s^z)^2}$, and $h$ denotes the distance between the structural center and the ocean surface.

In the shallow ocean, underwater acoustic propagation involves multiple reflections from the ocean surface and ocean sediment. To represent the reflection characteristic, the PIKF can be regarded as the superposition of multiple image point source acoustic waves, which can be described as

$$\Psi_3 = \sum_{\chi=0}^{\infty} (\beta_1\beta_2)^\chi \left( \frac{e^{-ikR_1}}{R_1} + \beta_1 \frac{e^{-ikR_2}}{R_2} + \beta_2 \frac{e^{-ikR_3}}{R_3} + \beta_1\beta_2 \frac{e^{-ikR_4}}{R_4} \right), \quad (8)$$

where

$$R_1 = \sqrt{(x-s^x)^2 + (y-s^y)^2 + (2\chi H + z - s^z)^2},$$
$$R_2 = \sqrt{(x-s^x)^2 + (y-s^y)^2 + (2\chi H + 2(H-h) + z + s^z)^2},$$
$$R_3 = \sqrt{(x-s^x)^2 + (y-s^y)^2 + (2\chi H + 2h - z - s^z)^2},$$
$$R_4 = \sqrt{(x-s^x)^2 + (y-s^y)^2 + (2(\chi+1)H - z + s^z)^2},$$

in which $H$ indicates the depth of the shallow ocean, $\beta_1$ and $\beta_2$ represent the acoustic reflection coefficients of the ocean sediment and ocean surface, respectively. The ocean surface is typically considered as a free surface, the acoustic reflection coefficient $\beta_2 = -1$. The acoustic reflection coefficient $\beta_1$ is determined by the mass density and sound speed of seawater and ocean sediment, which can be obtained by

$$\beta_1 = \frac{\rho_1 \cos(\theta)/\rho_0 - \sqrt{(c_0/c_1)^2 - \sin^2(\theta)}}{\rho_1 \cos(\theta)/\rho_0 + \sqrt{(c_0/c_1)^2 - \sin^2(\theta)}}, \quad (9)$$

where $\theta$ denotes the incident angle of the acoustic wave. $c_0$ and $\rho_0$ are the sound speed and mass density of the seawater, $c_1$ and $\rho_1$ represent the sound speed and mass density of the ocean sediment.

### 2.2.2. Numerical framework of PIKFNNs

Based on the PIKFs of underwater acoustic propagation in different ocean environments, this section presents the numerical framework of PIKFNNs. To provide a more intuitive presentation, we first introduce the well-known physics-informed neural networks (PINNs). The core idea of PINNs is to embed physical equations into loss functions, allowing the neural network to learn the underlying physical laws from the data, thereby improving model accuracy and generalization capabilities.

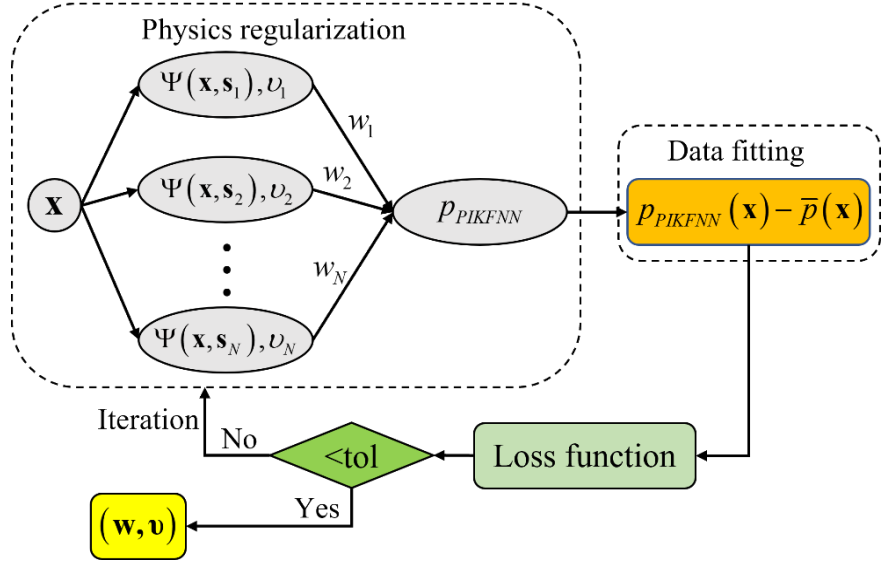

Fig. 2. Numerical framework of PIKFNNs.

Inspired by the above perspective, we replaced the conventional activation functions in neural networks with PIKFs, thus integrating prior physical information into the neural network and constructing the PIKFNNs, as illustrated in Fig. 2. The PIKFNNs is a fully connected neural architecture consisting of an input layer, a hidden layer, and an output layer. $\mathbf{w}$ and $\upsilon$ are the parameters in the neural network. $\Psi$ represents the PIKFs of underwater acoustic propagation in different ocean environments, which contains prior physical information and thus naturally performs the physical regularization in the hidden layer of the neural network.

For the underwater acoustic propagation in the ocean environment, the loss function of the PIKFNN can be written as

$$Loss_{PIKFNN} = \frac{1}{N}\sum_{i=1}^{N}\left(p_{PIKFNN}(\mathbf{x}_i) - \bar{p}(\mathbf{x}_i)\right)^2, \tag{10}$$

where $N$ is the number of training samples, $p_{PIKFNN}$ denotes the predicted sound pressure of PIKFNNs. Eq. (10) represents the mean squared error between the true solution and the predicted values of underwater sound pressure, which signifies the data fitting process in the neural network.

Then, the optimization algorithm is utilized to optimize the loss function and update the parameters of the neural network to gradually fit the boundary training data. In this study, the Levenberg-Marquardt (LM) algorithm is utilized as the optimizer for training the neural network model, which integrates the strengths of the Gauss-Newton method and the gradient descent method, and allowing automatic adjustment of the step size during iterations to balance convergence speed and stability. The optimization process of the LM algorithm is as follows.

First, the objective of the optimization algorithm is to minimize the loss function in PIKFNNs

$$\mathbf{u} = \arg\min_{u} Loss. \tag{11}$$

Then, the stopping tolerance $tol$ is determined, and the optimization result is checked to see if it is less than the tolerance

$$\textbf{if } \max|\mathbf{u}_i - \mathbf{u}_{i-1}| < tol \quad or \quad |Loss_i(\mathbf{u}_i) - Loss_i(\mathbf{u}_i)| < tol, \tag{12}$$

where the subscript '$i$' indicates the $i$-th iteration. Finally, the optimization is stopped when the result is less than the tolerance.

The above process constitutes the training process of PIKFNNs, and the trained PIKFNNs can be applied to new data to obtain predicted underwater sound pressure. It should be noted here that while both PINNs and PIKFNNs incorporate prior physical information into the neural network, PINNs embeds the governing equations into the loss function, while PIKFNNs replaces traditional activation functions with PIKFs as neuron basis functions, which is the main difference between the two neural network models.

## 3. Numerical examples and discussions

This section presents some examples to verify the effectiveness and accuracy of the developed method for predicting the underwater acoustic propagation induced by structural vibrations in the unbounded ocean, deep ocean and shallow ocean. The sound speed $c_0 = 1500\,m/s$ and mass density $\rho_0 = 1025\,kg/m^3$ of the seawater are considered. The L$_2$ relative error $Lrerr$ is defined by

$$Lreer = \sqrt{\sum_{i=1}^{N_t}(u_p(\mathbf{x}_i) - u_t(\mathbf{x}_i))^2 \Big/ \sum_{i=1}^{N_t}(u_t(\mathbf{x}_i))^2}, \tag{13}$$

where $N_t$ is the number of underwater test points, $u_p$ and $u_t$ indicate the predicted solution and true solution at the node $\mathbf{x}_i$, respectively.

**Example 1.** Unbounded ocean

To verify the effectiveness of the PIKFNN, we first study the underwater acoustic field induced by the vibration of a spherical structure in the unbounded ocean. The radius of the sphere is $r = 1m$, and the radial vibration velocity is $v_0 = 0.1 mm/s$, as shown in Fig. 3. The true solution for the underwater sound pressure in the unbounded ocean is

$$p(\bar{r}) = \frac{r}{\bar{r}}\left(\frac{ikr\rho_0 c_0}{ikr-1}\right)v_0 e^{ik(\bar{r}-r)}, \tag{14}$$

where $\bar{r}$ denotes the distance from the underwater test point to the center of the sphere, $k$ is the wave number.

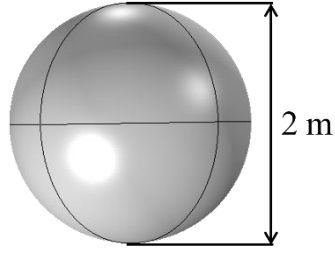

Fig. 3. Sketch of the spherical structure.

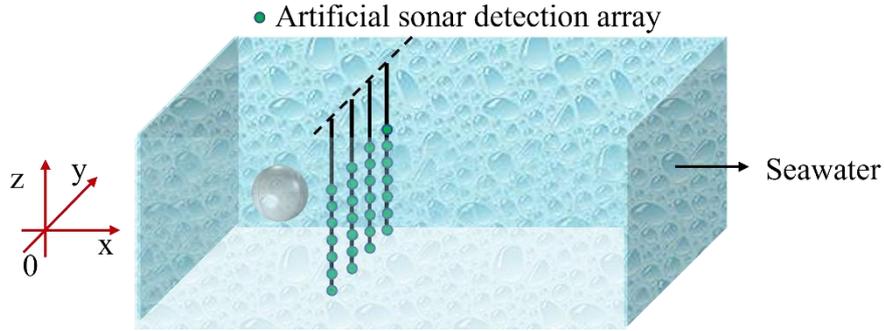

Fig. 4. Diagram of the artificial sonar detection array in the unbounded ocean.

In the PIKFNN implementation, the artificial sonar detection array is deployed at a distance of $3m$ from the center of the underwater structure, and the Cartesian coordinate is established with the center of the structure as the origin, as depicted in Fig. 4. The artificial sonar detection array consists of 9 linear array sonars distributed along the y-axis, with a spacing of $0.5m$ between the linear array sonars. Each linear array sonar is equipped with 17 hydrophones spaced $0.5m$ apart. The above true solution is considered to be the underwater sound pressure measured by the artificial sonar detection array. The PIKF function is $\Psi_1$, the source points are uniformly distributed on a sphere with a radius of $0.5m$, and the training data consists of underwater sound pressure obtained from the artificial sonar detection array. The number of PIKF neurons, source points, and training samples is $N=153$. The Levenberg-Marquardt algorithm is utilized to minimize the loss function formed by the training data. 1911 underwater test points are uniformly distributed in a rectangular area $\{(x,y,z)|10m \leq x \leq 100m, y=0, -10m \leq z \leq 10m\}$.

We investigate the influence of the tolerance and the number of training samples on the numerical accuracy of the proposed PIKFNN, the vibration frequency of the sphere is 6000Hz. Tables 1 and 2 compare the L$_2$ relative error of the proposed PIKFNN with different tolerances and number of training samples. It can be found that the PIKFNN results gradually converge as the tolerance decreases and the number of training samples increases. Here, the tolerance of $tol=1\text{E-}6$ and $N=153$ training samples are used for the subsequent PIKFNN computation.

Table 1. *Lrerr* with different tolerances.

| *tol* | 1E-1 | 1E-2 | 1E-3 | 1E-4 | 1E-5 | 1E-6 |
|---|---|---|---|---|---|---|
| *Lrerr* | 1.17E-4 | 1.02E-4 | 9.64E-5 | 6.64E-5 | 4.03E-5 | 3.61E-5 |

Table 2. *Lrerr* with different number of training samples.

| *N* | 85 | 119 | 153 | 187 | 221 |
|---|---|---|---|---|---|
| *Lrerr* | 5.78E-3 | 1.01E-3 | 3.61E-5 | 2.82E-5 | 1.08E-5 |

Fig. 5 illustrates the sound pressure levels predicted by the developed PIKFNN and the true solution at different vibration frequencies. The coordinates of the underwater test point are (10, 0, 0). The contour plots of the sound pressure level induced by the structural vibration frequency of 1000 Hz are shown in Fig. 6. It can be seen that the underwater sound pressure level results obtained by the developed PIKFNN fit well with the true solution. These figures further confirm the reliability of the developed PIKFNN for predicting the sound pressure levels induced by structural vibrations in the unbounded ocean.

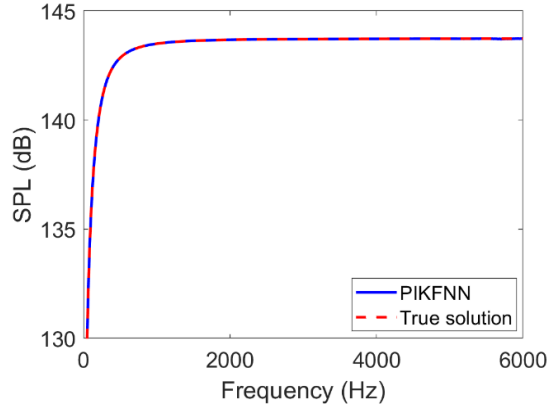

Fig. 5. The sound pressure levels obtained by the PIKFNN and true solution.

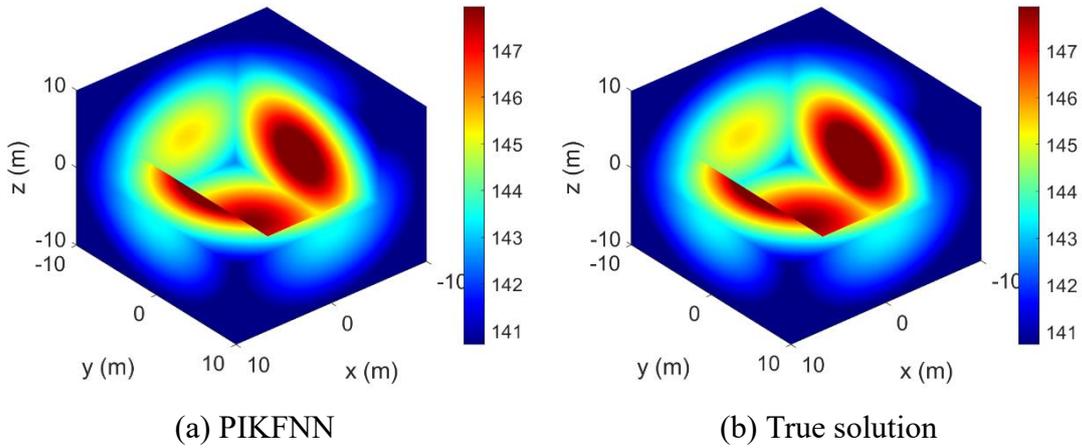

(a) PIKFNN  (b) True solution

Fig. 6. Contours of sound pressure level induced by the vibration of a spherical structure in the unbounded ocean.

**Example 2.** Deep ocean

As the second numerical example, the prediction performance of FEM-PIKFNN for the underwater acoustic propagation induced by the vibration of the ribbed cylindrical shell in the deep ocean environment is investigated. At the bottom of the ribbed cylindrical shell, there is a simple harmonic concentrated force $F_1 = 10N$ acting vertically downward, while at the middle and top, two simple harmonic distributed forces $F_2 = 10 N/m$ act outward along the horizontal direction. The specific dimensions and force diagram of the ribbed cylindrical shell are shown in Fig. 7. The immersion depth of the ribbed cylindrical shell is $h = 20m$.

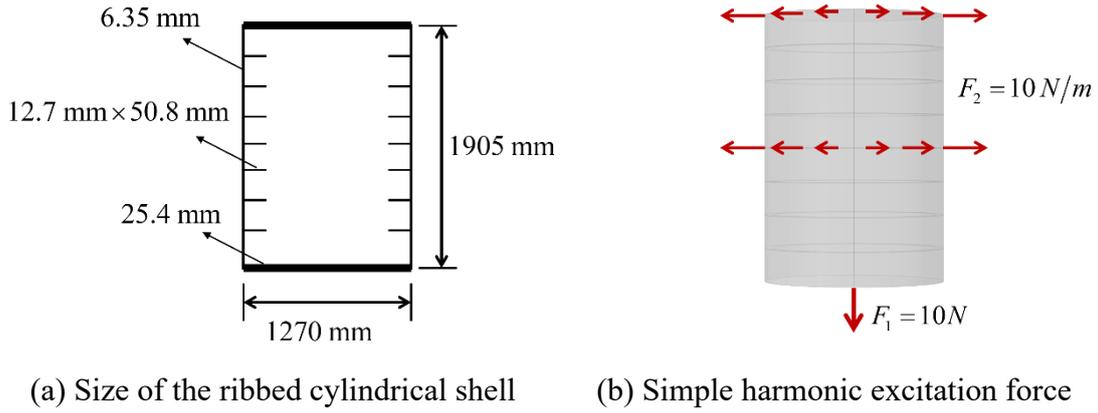

(a) Size of the ribbed cylindrical shell    (b) Simple harmonic excitation force

Fig. 7. Sketch of the ribbed cylindrical shell.

In the FEM-PIKFNN implementation, due to the axisymmetric characteristics of the structure shape and force, the FEM uses the axisymmetric model to calculate the near-field sound pressure level induced by the vibration of the ribbed cylindrical shell structure, with the number of triangular and quadrilateral elements being 143846 and 52800, respectively. The results calculated by the FEM are used as the underwater sound pressure obtained by the artificial sonar detection array. The position of the artificial sonar detection array is the same as in Example 1. The PIKF function is $\Psi_2$, 153 samples obtained from the artificial sonar detection array in the deep ocean are adopted for neural network training, and 153 source points are uniformly distributed on a sphere with a radius of $0.5m$. The Levenberg-Marquardt algorithm with a tolerance of $tol = 1\text{E-}6$ is used to minimize the loss function composed of underwater sound pressure in the deep ocean.

Fig. 8 plots the sound pressure levels at the underwater test nodes (10, 0, 0) calculated by the present FEM-PIKFNN and COMSOL Multiphysics finite element analysis (FEA) software. As illustrated in Fig. 8, the sound pressure levels obtained by both methods are consistent over the range of excitation frequencies from 50 Hz to 1000 Hz. Next, two different excitation frequencies, 400 Hz and 800 Hz, are considered. The contour plots of the sound pressure level induced by the vibration of the ribbed cylindrical shell in the deep ocean environment are plotted in Figs. 9 and 10, respectively. It can be observed that the sound pressure level results predicted by the present FEM-PIKFNN agree well with those calculated by the COMSOL Multiphysics FEA software for the simple harmonic force with two different excitation frequencies, which indicates that the present FEM-PIKFNN can achieve satisfactory results for predicting the underwater acoustic propagation in the deep ocean.

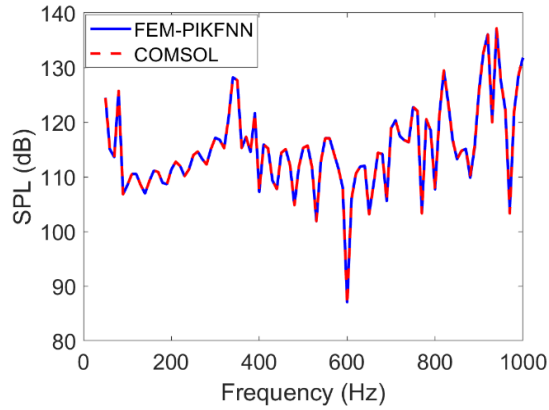

Fig. 8. The sound pressure levels in the deep ocean obtained by the FEM-PIKFNN and COMSOL Multiphysics.

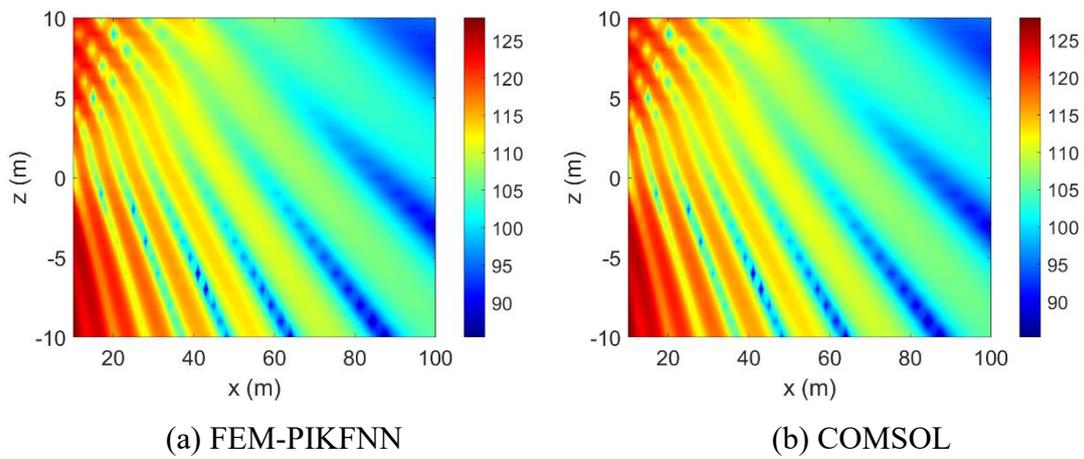

(a) FEM-PIKFNN    (b) COMSOL

Fig. 9. Contours of the sound pressure level of the ribbed cylindrical shell under simple harmonic force excitation with a frequency of 400Hz.

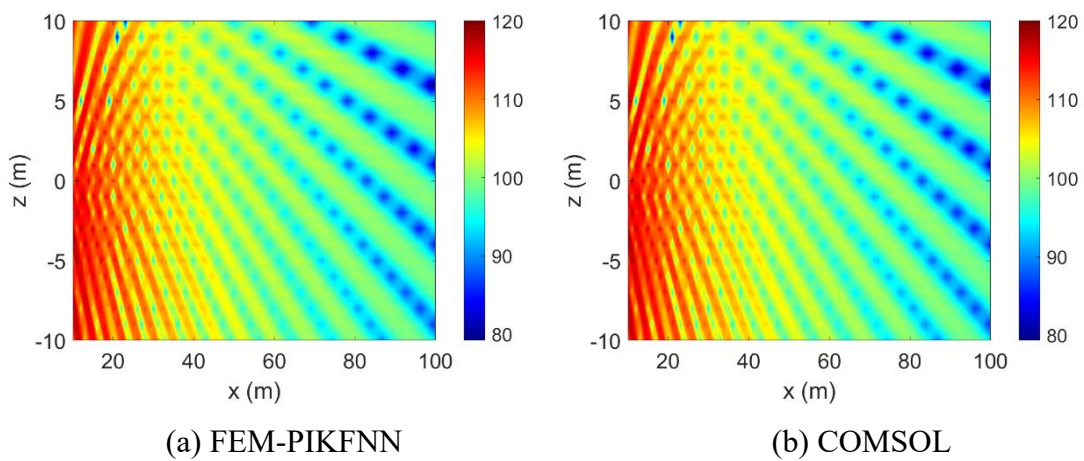

(a) FEM-PIKFNN    (b) COMSOL

Fig. 10. Contours of the sound pressure level of the ribbed cylindrical shell under simple harmonic force excitation with a frequency of 800Hz.

**Example 3.** Shallow ocean

In the last example, the underwater acoustic propagation behavior induced by the vibration of a capsule shell under simple harmonic force in the shallow ocean environment is considered. The two simple harmonic distributed force $F_3 = 5\,N/m$ acts outward along the horizontal direction on the capsule shell with thickness $th = 1mm$, and its sketch is presented in Fig. 11. The depth of the shallow ocean is $H = 20m$, and the immersion depth of the capsule shell is $h = 10m$. The mass density and sound speed of the ocean sediment are $\rho_1 = 2600\,kg/m^3$ and $c_1 = 1620\,m/s$, respectively.

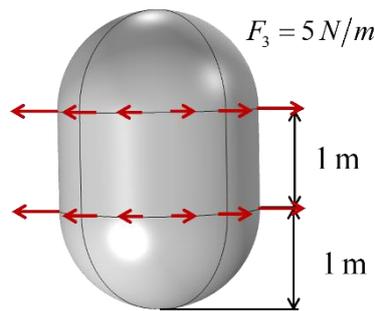

Fig. 11. Sketch of the capsule shell with simple harmonic excitation force.

Similar to Example 2, the underwater sound pressure induced by the vibration of the capsule shell structure in the shallow ocean calculated by the FEM is used as the training data of PIKFNN. The number of triangular and quadrilateral elements in the finite element simulation is 34048 and 22400, respectively. In the PIKFNN implementation, the PIKF function is $\Psi_3$, 153 training data are obtained from the artificial sonar detection array in the shallow ocean, and 153 source points are uniformly distributed on a sphere with a radius of $0.5m$. The Levenberg-Marquardt algorithm with a tolerance of $tol = 1E\text{-}6$ is used to minimize the loss function in this example.

The underwater sound pressure level curve induced by the capsule shell at different excitation frequencies in the shallow ocean environment is depicted in Fig. 12. It can be found that the sound pressure levels predicted by the FEM-PIKFNN are in good agreement with those calculated by the COMSOL Multiphysics FEA software at the underwater test nodes (10, 0, 0). The underwater sound pressure level distributions of the capsule shell excited by two distributed forces with the frequency of 500 Hz at the immersion depths of $h = 8m$ and $h = 12m$ are depicted in Figs. 13 and 14, respectively. The numerical results demonstrate that the predicted sound pressure level obtained by the proposed FEM-PIKFNN is consistent with that calculated by the COMSOL Multiphysics FEA software, which shows that the FEM-PIKFNN is effective in predicting the underwater acoustic propagation behavior in the shallow ocean environment.

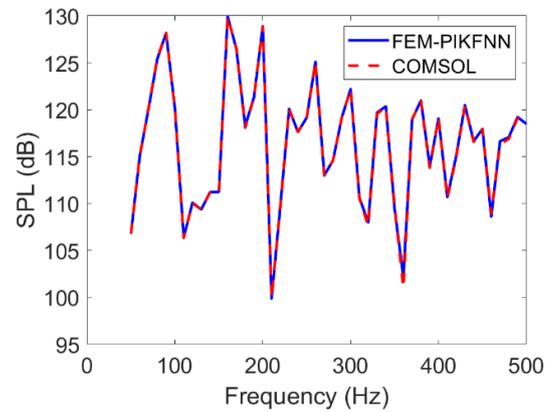

Fig. 12. The sound pressure levels in the shallow ocean obtained by the FEM-PIKFNN and COMSOL Multiphysics.

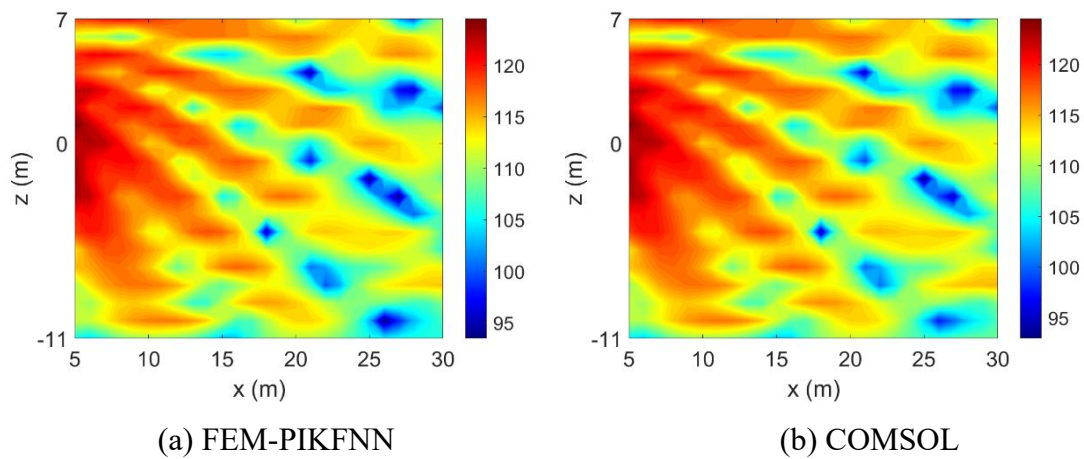

(a) FEM-PIKFNN  (b) COMSOL

Fig. 13. Contours of the sound pressure level of the capsule shell ($h = 8m$).

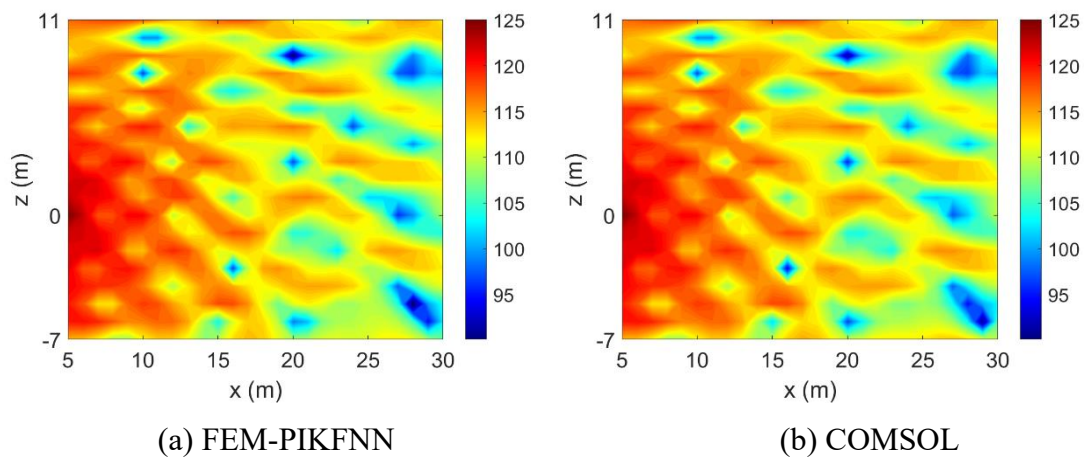

(a) FEM-PIKFNN  (b) COMSOL

Fig. 14. Contours of the sound pressure level of the capsule shell ($h = 12m$).

## 5. Conclusions

This paper proposes the hybrid FEM-PIKFNNs for predicting the underwater acoustic propagation induced by structural vibrations in the unbounded ocean, deep ocean and shallow ocean. In the proposed method, the underwater sound pressure levels calculated by FEM are used as the training data of PIKFNNs. PIKFs containing the underwater acoustic propagation information from different ocean environments are used as the activation functions of PIKFNNs, and the loss function contains only the training sample information of the underwater acoustic propagation. The numerical results show that PIKFNNs gradually converge as the tolerance decreases and the number of training samples increases, and the effectiveness and accuracy of the proposed method in predicting underwater acoustic propagation are verified by comparison with the true solutions and finite element results.

The study demonstrates that incorporating prior physical information into the activation functions of artificial neural networks can improve their generalization ability and reduce the complexity of neural networks. PIKFNNs with a single hidden layer can accurately predict the underwater acoustic propagation behavior in three different ocean environments.

**Acknowledgements**
The work reported in this paper was supported by the Natural Science Foundation of China (Grant No. 12122205), the Innovative Research Foundation of Ship General Performance (Grant No. 33122126), the Six Talent Peaks Project in Jiangsu Province of China (Grant No. 2019-KTHY-009).